\documentclass{amsart}

\usepackage{amsmath,amssymb,amscd,amsxtra}
\usepackage{latexsym}
\usepackage[dvips]{graphics,epsfig}

\newtheorem{theorem}{Theorem}[section]
\newtheorem{lemma}[theorem]{Lemma}
\newtheorem{proposition}[theorem]{Proposition}
\newtheorem{cor}[theorem]{Corollary}

\theoremstyle{definition}
\newtheorem{definition}[theorem]{Definition}
\newtheorem{example}[theorem]{Example}

\theoremstyle{remark}
\newtheorem{remark}[theorem]{Remark}

\numberwithin{equation}{section}

\newcommand{\Bc}{\mathbb C}
\newcommand{\Bn}{\mathbb N}
\newcommand{\Br}{\mathbb R}
\newcommand{\Bt}{\mathbb T}
\newcommand{\Bz}{\mathbb Z}
\newcommand{\Bh}{\mathbb H}

\allowdisplaybreaks

\begin{document}

\title{Frame potential and finite abelian groups}

\author{Brody D. Johnson}
\address{Department of Mathematics and Computer Science,
Saint Louis University, Saint Louis, MO 63103}
\email{brody@slu.edu}

\author{Kasso A. Okoudjou}
\address{Department of Mathematics,
University of Maryland, College Park, MD 20742}
\email{kasso@math.umd.edu}

\date{\today}

\begin{abstract}
This article continues a prior investigation of the authors with the goal of extending characterization results of convolutional tight frames from the context of cyclic groups to general finite abelian groups.  The collections studied are formed by translating a number of \emph{generators} by elements of a fixed subgroup and it is shown, under certain norm conditions, that tight frames with this structure are characterized as local minimizers of the frame potential.  Natural analogs to the downsampling and upsampling operators of finite cyclic groups are studied for arbitrary subgroups of finite abelian groups.  Directions of further study are also proposed.
\end{abstract}

\maketitle

\section{Introduction} \label{section1}

In recent years, there has been much activity in the study of frames for finite-dimensional Hilbert spaces.  Topics of interest include the characterization and construction of frames of various kinds, e.g., equiangular frames \cite{Renes}, harmonic frames \cite{ValeWaldron2}, [compound] geometrically uniform frames \cite{BE}, tight frames \cite{BF,CFKLT}, frames with symmetries \cite{ValeWaldron1}, frames resistant to erasures \cite{HolmesPaulsen}, etc.  Throughout this section, let $\Bh$ denote a finite dimensional real or complex Hilbert space.  Recall that a \emph{frame} for $\Bh$ is a collection $X=\lbrace f_{k} \rbrace_{k=1}^{n} \subseteq \Bh$ for which there exist real numbers $0<A\le B<\infty$ such that
\begin{equation*}
A \Vert f\Vert^{2} \le \sum_{k=1}^{n} \vert \langle f, f_{k} \rangle \vert^{2} \le B \Vert f \Vert^{2}, \quad \text{for all} \; f\in \Bh.
\end{equation*}

\noindent
If it is possible to choose $A=B$ then $X$ is called a \emph{tight frame}.

One recent line of research began with the work of Benedetto and Fickus in \cite{BF}, where a characterization of tight frames composed of unit-norm vectors was given in terms of the frame potential.

\begin{definition} \label{FPdef}
Let $\Bh$ be a finite dimensional real or complex Hilbert space.  Let $X=\lbrace f_{k} \rbrace_{k=1}^{N} \subseteq \Bh$, then the \emph{frame potential} of the collection $X$ is defined as
\begin{equation*}
\text{FP}(X) = \sum_{j=1}^{N} \sum_{k=1}^{N} \left \vert \left \langle f_{j}, f_{k} \right \rangle \right \vert^{2}.
\end{equation*}
\end{definition}

Following the work of Benedetto and Fickus, Casazza et al. \cite{CFKLT} considered frames composed of vectors with arbitrary norms and arrived at the following description of tight frames.

\begin{theorem}[Theorem 10 of \cite{CFKLT}] \label{FFthm}
Let $\Bh$ be a $d$-dimensional real or complex Hilbert space and fix $a_{0}\ge a_{1} \ge \cdots \ge a_{M-1}>0$, $M\ge d$.  Denote by $m_{0}$ the smallest index $0\le m\le M-1$ for which
\begin{equation} \label{ineq}
(d-m) a_{0}^{2} \le \sum_{j=m}^{M-1} a_{j}^{2}
\end{equation}

\noindent
holds.  If $X=\lbrace f_{j} \rbrace_{j=0}^{M-1} \subset \Bh$ is a local minimizer of the frame potential subject to the constraint $\Vert f_{j} \Vert = a_{j}$, $0\le j\le M-1$, then $X$ may be divided into two mutually orthogonal subcollections: $\lbrace f_{j} \rbrace_{j=0}^{m_{0}-1}$, which consists of mutually orthogonal nonzero vectors, and $\lbrace f_{j} \rbrace_{j=m_{0}}^{M-1}$, which is a tight frame for its $(d-m_{0})$-dimensional span.  In the event that $m_{0}=0$, $X$ is a tight frame for $\Bh$.
\end{theorem}

When $m=0$, \eqref{ineq} is referred to as the \emph{fundamental frame inequality},
\begin{equation} \label{FFineq}
d a_{0}^{2} \le \sum_{j=0}^{M-1} a_{j}^{2}.
\end{equation}

\noindent
In light of Theorem \ref{FFthm}, it is clear that this inequality provides a sufficient condition for the existence of a tight-frame having a specified nonincreasing sequence of norms.

\begin{remark}
The following additional facts from \cite{CFKLT} related to Theorem \ref{FFthm} and the inequality \eqref{FFineq} will also be relevant to this work.
\renewcommand{\labelenumi}{(\roman{enumi})}
\begin{enumerate}
\item  The fundamental frame inequality \eqref{FFineq} is a necessary condition on the norms associated to any tight frame after rearrangement into decreasing order.

\item  Any local minimizer of the frame potential (associated with a fixed sequence of norms) must also be a global minimizer.
\end{enumerate}
\end{remark}

Another line of research, originated by Vale and Waldron \cite{ValeWaldron1}, deals with an examination of certain symmetries possessed by tight frames for finite dimensional Hilbert spaces.  Following \cite{ValeWaldron1}, define the \emph{symmetry group} of a frame $X=\lbrace f_{k} \rbrace_{k=1}^{n}$ for $\Bh$ to be the group
\begin{equation*}
\text{Sym}(X) = \lbrace U \in \mathcal{U}(\Bh) : U(X) = X \rbrace.
\end{equation*}

\noindent
Here, $\mathcal{U}(\Bh)$ denotes the group of unitary linear transformations on $\Bh$ under composition.  This motivates a natural question: \emph{Under what conditions, if any, can tight frames with specified symmetries be characterized as local minimizers of the frame potential?}

\begin{example}
Let $\omega = e^{2\pi i/3}$ and define $A$ to be the $2\times 2$ matrix
\begin{equation*}
A = \begin{bmatrix} \omega & 0 \\ 0 & \omega^{2} \end{bmatrix}.
\end{equation*}

\noindent
Given $u\in \Bc^{2}$ with $\Vert u\Vert=1$, consider the collection $X=\lbrace A^{j}u\rbrace_{j=0}^{2}$.  It will be shown that the local minimizers of the frame potential of $X$ (under the constraint that $\Vert u\Vert =1$) are precisely the tight frames of this form.  By symmetry,
\begin{equation*}
\text{FP}(X) = 3 \sum_{k=0}^{2} \left \vert \left \langle A^{k} u, u \right \rangle \right \vert^{2},
\end{equation*}

\noindent
so letting $u=(u_{1},u_{2})$ one seeks to minimize (after an elementary computation)
\begin{equation*}
f(u_{1},u_{2}) = 3 \left (1 + 2u_{1}^{4} - 2u_{1}^{2} u_{2}^{2} + 2u_{2}^{4} \right )
\end{equation*}

\noindent
subject to the constraint that $g(u_{1},u_{2}) = u_{1}^{2} + u_{2}^{2} = 1$.  The method of Lagrange multipliers reveals that the minima occur when $u_{1}^{2}=u_{2}^{2}$, i.e., $u$ is of the form $(\pm 1/\sqrt{2}, \pm 1/\sqrt{2})$ or $(\pm 1/\sqrt{2}, \mp 1/\sqrt{2})$.  Each of these four choices leads to a $\frac{3}{2}$-tight frame for $\Bc^{2}$.  Hence, each local minimizer of $\text{FP}(X)$ leads to a tight frame.  Moreover, if there were another choice for $u$ which led to a tight frame it would achieve the same minimum value for $f$ and, therefore, would be among the solutions found through Lagrange multipliers.  This shows that the tight frames of the form $\lbrace A^{j}u\rbrace_{j=0}^{2}$ are precisely the local minimizers of the frame potential.
\end{example}

The present work seeks an answer to the preceding question for certain symmetry groups in association with the group algebra $\ell(G)$, the real or complex Hilbert space of functions defined on a finite abelian group $G$.  Notice that each element $g\in G$ leads to a natural \emph{translation} operator on $\ell(G)$,
\begin{equation*}
T_{g}: \ell(G) \rightarrow \ell(G), \quad (T_{g} f)(g') = f(g'g^{-1}).
\end{equation*}

\noindent
Therefore, given a subgroup $H$ of $G$, any frame of the form
\begin{equation} \label{convsystem}
X_{H} = \lbrace T_{h}f_{k} : h\in H, \, f_{k}\in \ell(G),\, 0\le k\le n-1 \rbrace
\end{equation}

\noindent
satisfies $H \le \text{Sym}(X_{H})$.  Moreover, such frames bear a close relationship with the convolution structure of $\ell(G)$, which leads to efficient implementation in applications by way of the fast Fourier transform.  In this sense, this work continues an investigation of the authors' from \cite{FJKO} in which tight frames generated by translations in $\ell(\Bz/d\Bz)$ were characterized as the local minimizers of the frame potential under certain norm conditions related to \eqref{ineq}.

\section{Preliminaries} \label{section2}

Throughout this section $G$ will denote a finite abelian group.  Recall that the inner product on $\ell(G)$ is given by
\begin{equation*}
\langle f_{1}, f_{2} \rangle = \sum_{g\in G} f_{1}(g) \overline{f_{2}(g)}, \quad f_{1},f_{2}\in \ell(G).
\end{equation*}

\noindent
The convolution of $f_{1}, f_{2} \in \ell(G)$, denoted $f_{1}*f_{2}\in \ell(G)$, is given by
\begin{equation*}
f_{1}*f_{2}(g) = \sum_{x\in G} f_{1}(x)f_{2}(g^{-1}x), \quad g\in G.
\end{equation*}

This work is concerned with the study of collections of \emph{filters}, $X = \lbrace f_{m} \rbrace_{m=0}^{n-1} \in \ell(G)$, which will be used to analyze and synthesize general elements of $\ell(G)$ via convolution.  It is natural to \emph{sample} such convolutions over a subgroup $H$, which may be loosely interpreted as a downsampling operation with respect to the quotient group $G/H$.  With this motivation, define the \emph{sampling operator over $H$}, $\mathcal{S}_{H} : \ell(G) \rightarrow \ell(H)$, by
\begin{equation*}
(\mathcal{S}_{H}f)(h) = f(h), \quad h\in H.
\end{equation*}

\noindent
Similarly, define the \emph{upsampling operator over $H$}, $\mathcal{S}^{*}_{H} : \ell(H)\rightarrow \ell(G)$, by
\begin{equation*}
(\mathcal{S}^{*}_{H}f)(g) = \begin{cases} f(g), & g\in H \\ 0, & g\notin H, \end{cases} \quad g\in G.
\end{equation*}

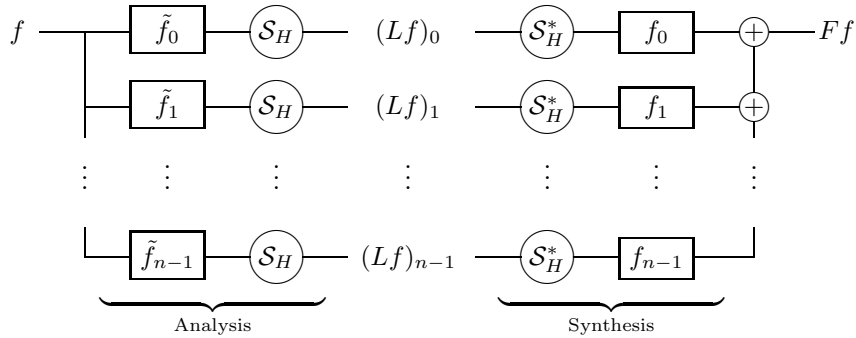
\begin{figure}[hbtp]
\begin{center}
  \begin{picture}(110,42)(0,0)
    \put (1,39){$f$}
    \put (5,40){\line(1,0){12}}
    \put (11,40){\line(0,-1){14}}
    \put (11,10){\line(0,1){4}}

    \put (17,39){\framebox[10mm]{$\tilde{f}_{0}$}}
    \put (27,40){\line(1,0){6}}
    \put (36.5,40){\circle{7}}
    \put (33,39){\makebox[7mm]{$\mathcal{S}_{H}$}}
    \put (40,40){\line(1,0){6}}
    \put (49,39){\makebox[10mm]{$(Lf)_{0}$}}
    \put (63,40){\line(1,0){6}}
    \put (72.5,40){\circle{7}}
    \put (69,39){\makebox[7mm]{$\mathcal{S}^{*}_{H}$}}
    \put (76,40){\line(1,0){6}}
    \put (82,39){\framebox[10mm]{$f_{0}$}}

    \put (11,30){\line(1,0){6}}
    \put (17,29){\framebox[10mm]{$\tilde{f}_{1}$}}
    \put (27,30){\line(1,0){6}}
    \put (36.5,30){\circle{7}}
    \put (33,29){\makebox[7mm]{$\mathcal{S}_{H}$}}
    \put (40,30){\line(1,0){6}}
    \put (49,29){\makebox[10mm]{$(Lf)_{1}$}}
    \put (63,30){\line(1,0){6}}
    \put (72.5,30){\circle{7}}
    \put (69,29){\makebox[7mm]{$\mathcal{S}^{*}_{H}$}}
    \put (76,30){\line(1,0){6}}
    \put (82,29){\framebox[10mm]{$f_{1}$}}

    \put (10,19){\makebox[2mm]{$\vdots$}}
    \put (17,19){\makebox[10mm]{$\vdots$}}
    \put (33,19){\makebox[7mm]{$\vdots$}}
    \put (49,19){\makebox[10mm]{$\vdots$}}
    \put (69,19){\makebox[7mm]{$\vdots$}}
    \put (82,19){\makebox[10mm]{$\vdots$}}

    \put (11,10){\line(1,0){6}}
    \put (17,9){\framebox[10mm]{$\tilde{f}_{n-1}$}}
    \put (27,10){\line(1,0){6}}
    \put (36.5,10){\circle{7}}
    \put (33,9){\makebox[7mm]{$\mathcal{S}_{H}$}}
    \put (40,10){\line(1,0){6}}
    \put (49,9){\makebox[10mm]{$(Lf)_{n-1}$}}
    \put (63,10){\line(1,0){6}}
    \put (72.5,10){\circle{7}}
    \put (69,9){\makebox[7mm]{$\mathcal{S}^{*}_{H}$}}
    \put (76,10){\line(1,0){6}}
    \put (82,9){\framebox[10mm]{$f_{n-1}$}}

    \put (92,40){\line(1,0){6}}
    \put (100,40){\circle{4}}
    \put (98,39){\makebox[4mm]{$+$}}
    \put (102,40){\line(1,0){6}}
    \put (108,39){\makebox[6mm]{$Ff$}}
    \put (100,32){\line(0,1){6}}
    \put (92,30){\line(1,0){6}}
    \put (100,26){\line(0,1){2}}
    \put (100,30){\circle{4}}
    \put (98,29){\makebox[4mm]{$+$}}
    \put (98,19){\makebox[4mm]{$\vdots$}}
    \put (100,10){\line(0,1){4}}
    \put (92,10){\line(1,0){8}}

    \put (13,5){$\underset{\text{Analysis}}{\underbrace{\makebox[30mm]{}}}$}
    \put (66,5){$\underset{\text{Synthesis}}{\underbrace{\makebox[30mm]{}}}$}

  \end{picture}
\end{center}
\caption{Block diagram of an $n$-channel filterbank on $\ell(G)$.} \label{filterbank}
\end{figure}

\noindent
Figure \ref{filterbank} depicts a typical filterbank, composed of analysis and synthesis stages that make use of convolution (represented by rectangular elements) as well as the sampling and upsampling operators (represented by circular elements).  The analysis stage implements the \emph{involution} of each filter, an operation on $\ell(G)$ given by $\tilde{f}(g) = \overline{f(g^{-1})}$.

Associated with the filterbank of Figure \ref{filterbank} are several important operators.  The \emph{filterbank analysis operator}, $L:\ell(G)\rightarrow \bigoplus_{m=0}^{n-1} \ell(H)$, is given by
\begin{equation*}
Lf = (Lf)_{0} \oplus \cdots \oplus (Lf)_{n-1},
\end{equation*}

\noindent
where $(Lf)_{m} = \mathcal{S}_{H}(f*\tilde{f}_{m})$.  The adjoint of the filterbank analysis operator, $L^{*}:\oplus_{m=0}^{n-1} \ell(H)\rightarrow \ell(G)$, is called the \emph{filterbank synthesis operator} and is described by
\begin{equation*}
L^{*}\left ( \oplus_{m=0}^{n-1} y_{m} \right ) = \sum_{m=0}^{n-1} (\mathcal{S}^{*}_{H} y_{m})*f_{m}.
\end{equation*}

\noindent
The composition of the filterbank synthesis and analysis operators, $L^{*} L$, is called the \emph{filterbank frame operator} for the filterbank and is denoted by $F: \ell(G) \rightarrow \ell(G)$.  The frame operator is described explicitly by
\begin{equation*}
Ff = \sum_{m=0}^{n-1} \left \lbrack \mathcal{S}^{*}_{H} \mathcal{S}_{H}(f*\tilde{f}_{m}) \right \rbrack* f_{m}.
\end{equation*}

\noindent
The reader should note that the filterbank frame operator associated with $X$ is precisely the ordinary frame operator associated with the collection $X_{H}$, as defined by \eqref{convsystem}.  Of particular interest are collections of filters for which the associated filterbank frame operator is a scalar multiple of the identity, i.e., filters which give rise to tight frames for $\ell(G)$.  With this in mind, the main goal of this work is to provide a characterization of filterbanks which give rise to tight frames in terms of the frame potential.  Although neither of the works \cite{BF,CFKLT} includes the use of filterbanks, one could interpret their results in terms of a filterbank for $\ell(G)$ with sampling over the trivial subgroup.  (In this case, the group structure of $G$ plays no role in the analysis.)  As noted in Section \ref{section1}, the results of \cite{FJKO} correspond to a non-trivial filterbank with $G=\Bz_{d}$ and $H$ a cyclic subgroup of $G$.

The methods of \cite{FJKO} rely on a block diagonalization of $L^{*}$ called the \emph{modulated filter representation}, from which it follows that the frame characteristics of a filterbank for $\ell(\Bz_{d})$ are equivalent to the combined frame characteristics of a series of $d/N$ trivial filterbanks (downsampled over the trivial subgroup) for $\ell(\Bz_{N})$.  The characterization of tight frames for the latter case is described by the following theorem.

\begin{theorem}[Theorem 12 of \cite{FJKO}] \label{FJKOtheorem12}
Let $\lbrace a_{m} \rbrace_{m=0}^{n-1} \subset \Br$ be such that $a_{0} \ge a_{1} \ge \cdots \ge a_{n-1}>0$.  Let $d$ and $N$ be positive integers such that $N\mid d$ and $N\le n$. Denote by $m_{0}$ the smallest index $0\le m \le N-1$ such that
\begin{equation}  \label{cosetFFineq}
(N-m)a_{m}^{2} \le \sum_{j=m}^{n-1} a_{j}^{2}.
\end{equation}

\noindent
If the collections $Y_{j}:=\lbrace y_{m,j} \rbrace_{m=0}^{n-1} \subset \ell(\Bz_{N})$ form a local minimizer of the combined frame potential, $\sum_{j=0}^{\frac{d}{N}-1} \mathrm{FP}(Y_{j})$, under the constraint that
\begin{displaymath}
\sum_{j=0}^{\frac{d}{N}-1} \Vert y_{m,j} \Vert^{2} = \frac{d}{N} a_{m}^{2}, \quad 0\le m \le n-1,
\end{displaymath}

\noindent
then each collection $Y_{j}$ may be divided into two mutually orthogonal subcollections of $\ell(\Bz_{N})$: $\lbrace y_{m,j}\rbrace_{m=0}^{m_{0}-1}$, which consists of mutually orthogonal, nonzero vectors, and $\lbrace y_{m,j}\rbrace_{m=m_{0}}^{n-1}$, which is a tight frame for its $(N-m_{0})$-dimensional span.  Moreover, for each $j$ the norms of the vectors of $Y_{j}$ must satisfy $\Vert y_{m,j}\Vert=a_{m}$ for $0\le m \le m_{0}-1$ and $\sum_{m=m_{0}}^{n-1} \Vert y_{m,j}\Vert^{2} = \sum_{m=m_{0}}^{n-1} a_{m}^{2}$.  In the event that $m_{0}=0$ each collection $Y_{j}$, $0\le j \le \frac{d}{N}-1$, is a tight frame for $\ell(\Bz_{N})$ with a common frame bound.
\end{theorem}

The reader should note the similarity of \eqref{cosetFFineq} to \eqref{ineq}.  It is apparent from the statement of Theorem \ref{FJKOtheorem12} that this result is somewhat technical, yet it will play an important role in the analysis of filterbanks with sampling over arbitrary subgroups $H$.  Before such an analysis can be made, however, the behavior of $\mathcal{S}_{H}$ and $\mathcal{S}_{H}^{*}$ under the discrete Fourier transform will be examined.

\section{The discrete Fourier transform} \label{section3}

Throughout this section $G$ will represent a finite abelian group.  Recall that a \emph{character} of $G$ is a group homomorphism $\chi: \ell(G) \rightarrow \Bt$, where $\Bt$ represents the multiplicative group of unimodular complex numbers.  The \emph{dual group} to $G$ is denoted by $\widehat{G}$ and consists of all characters of $G$ under pointwise multiplication.  By the duality theorem of Pontryagin $\widehat{G}$ is, in fact, isomorphic to $G$, a result which can be obtained here as an easy consequence of the fundamental theorem of finite abelian groups.

Indeed, $G$ is isomorphic to a direct sum of cyclic groups, i.e.,
\begin{equation*}
G \simeq (\Bz/m_{1}\Bz) \oplus \cdots \oplus (\Bz/m_{r}\Bz),
\end{equation*}

\noindent
with $m_{j}\mid m_{j+1}$, $1\le j\le r-1$.  To each $a=(a_{1},\ldots,a_{r})\in G$ there is an associated character, $\chi_{a}$, given by
\begin{equation} \label{chardef}
\chi_{a}(x) = \prod_{j=1}^{r} \chi_{a_{j}}(x_{j}), \quad x=(x_{1},\ldots,x_{r}) \in G,
\end{equation}

\noindent
where $\chi_{a_{j}}(x_{j}) = \exp{(2\pi i a_{j}x_{j}/m_{j})}$.  For $a,b\in G$, therefore, it is evident that $\chi_{a} \chi_{b} = \chi_{ab}$.  The fact that the characters given by \eqref{chardef} exhaust $\widehat{G}$ is a consequence of the following lemma.

\begin{lemma}[\cite{Terras}] \label{ONcharacters}
Let $\chi, \psi \in \widehat{G}$, then
\begin{equation*}
\langle \chi, \psi \rangle = \begin{cases} \vert G \vert, & \chi=\psi, \\ 0, & \text{otherwise.} \end{cases}
\end{equation*}
\end{lemma}

The preceding discussion illustrates the fact that the characters of $G$ form an orthogonal basis of $\ell(G)$ and, consequently, that any $f\in \ell(G)$ is uniquely determined by its inner products with the characters.  This notion is the foundation for the discrete Fourier transform on $G$.

\begin{definition} \label{DFTabelian}
The \emph{discrete Fourier transform} (DFT) of $f\in \ell(G)$ is defined by
\begin{equation*}
\mathcal{F} f(\chi) = \hat{f}(\chi) = \sum_{x\in G} f(x)\overline{\chi (x)}, \quad \chi \in \widehat{G}.
\end{equation*}
\end{definition}

The following lemma summarizes the basic properties of the DFT on a finite abelian group, $G$.

\begin{lemma}[\cite{Terras}] \label{DFTproperties}
Basic Properties of the DFT.
\renewcommand{\labelenumi}{(\roman{enumi})}
\begin{enumerate}
\item $\mathcal{F}: \ell(G)\rightarrow \ell(\widehat{G})$ is a bijective linear map.

\item For $f_{1},f_{2}\in \ell(G)$, the convolution,
\begin{equation*}
(f_{1}*f_{2})(x) = \sum_{y\in G} f_{1}(y) f_{2}(x-y), \quad x\in G,
\end{equation*}

\noindent
satisfies
\begin{equation*}
\mathcal{F}(f_{1}*f_{2})(\chi) = \mathcal{F}f_{1}(\chi) \mathcal{F}f_{2}(\chi), \quad \chi \in \widehat{G}.
\end{equation*}

\item For $f\in \ell(G)$,
\begin{equation*}
f(x) = \frac{1}{\vert G \vert} \sum_{\chi \in \widehat{G}} \mathcal{F}f(\chi) \chi(x), \quad x\in G.
\end{equation*}

\item For $f_{1},f_{2}\in \ell(G)$,
\begin{equation*}
\langle f_{1}, f_{2} \rangle = \frac{1}{\vert G\vert} \langle \mathcal{F}f_{1}, \mathcal{F}f_{2} \rangle.
\end{equation*}
\end{enumerate}
\end{lemma}

Another important aspect of Fourier analysis on finite abelian groups, particularly for  the problems considered in this work, lies in the relationship between the dual groups of $G$ and a subgroup $H\le G$.  The following proposition is adapted from Proposition 6.1 of \cite{Serre}.

\begin{proposition} \label{Hext}
Suppose $H\le G$, let $x\in G\setminus H$, and denote by $H^{x}$ the subgroup of $G$ generated by $H$ and $x$.  Let $m_{x}=\min{\lbrace n\in \Bn : x^{n}\in H\rbrace}$.  Then each $\chi \in \widehat{H}$ extends to $m_{x}$ orthogonal characters in $\widehat{H^{x}}$, $\lbrace \chi_{j} \rbrace_{j=0}^{m_{x}-1}$, and
\begin{equation*}
\widehat{H^{x}} = \lbrace \chi_{j} : \chi_{j} \vert_{H} = \chi, \chi \in \widehat{H}, \; 0\le j \le m_{x}-1 \rbrace.
\end{equation*}
\end{proposition}

\begin{proof}
Fix $x\in G\setminus H$ and let $m:=m_{x}$.  Observe that $H^{x}$ is given by
\begin{equation*}
H^{x} = \lbrace x^{k}h : 0\le k \le m-1, \; h\in H\rbrace.
\end{equation*}

\noindent
If $x^{k}h_{1} = x^{\ell}h_{2}$ for some $0\le k,\ell < m$ and $h_{1},h_{2}\in H$ (without loss, assume $k\ge \ell$), then $x^{k-\ell} \in H$ with $0\le k-\ell < m$, a contradiction.  It follows that $\lbrack H^{x} : H\rbrack = m$.

Fix $\chi\in \widehat{H}$ and let $\omega = \chi (x^{m})$.  Define $\alpha_{j}$, $0\le j\le m-1$, as the distinct solutions of
\begin{equation*}
\alpha_{j}^{m} = \omega.
\end{equation*}

\noindent
Note that the collection $\lbrace \alpha_{j} \rbrace_{j=0}^{m-1}$ consists of a constant multiple of the $m$th roots of unity.  Now define $\chi_{j}$, $0\le j\le m-1$, by $\chi_{j}(x^{k}h) = \alpha_{j}^{k} \chi(h)$.  It is routine to verify that these elements of $\ell(H^{x})$ are characters.  Moreover, it follows that
\begin{equation*}
\langle \chi_{j_{1}}, \chi_{j_{2}} \rangle_{H^{x}} = \sum_{h\in H} \vert \chi(h)\vert^{2} \sum_{k=0}^{m-1} \alpha_{j_{1}}^{k} \, \overline{\alpha_{j_{2}}^{k}}.
\end{equation*}

\noindent
Hence, the orthogonality of $\lbrace \chi_{j} \rbrace_{j=0}^{m-1}$ follows from the orthogonality of the characters of $\Bz/m\Bz$.  Thus, each character of $\widehat{H}$ extends to $m$ orthogonal characters of $\widehat{H^{x}}$.

If $\chi,\psi\in \widehat{H}$ then
\begin{equation*}
\langle \chi_{j_{1}}, \psi_{j_{2}} \rangle_{H^{x}} = \sum_{h\in H} \chi(h) \, \overline{\psi(h)} \sum_{k=0}^{m-1} \alpha_{j_{1}}^{k} \, \overline{\tilde{\alpha}_{j_{2}}^{k}}
\end{equation*}

\noindent
and the orthogonality follows from the orthogonality of characters in $\widehat{H}$.  Dimensional considerations reveal that $\widehat{H^{x}}$ consists precisely of the character extensions claimed in the statement of the proposition.
\end{proof}

Let $H\le G$.  Given $\chi \in \widehat{H}$, let $\widehat{G}_{\chi}$ consist of all the characters in $\widehat{G}$ whose restrictions to $H$ coincide with $\chi$, i.e.,
\begin{equation*}
\widehat{G}_{\chi} = \lbrace \psi \in \widehat{G} : \psi\vert_{H} = \chi \rbrace.
\end{equation*}

\begin{cor} \label{Gchi}
Let $H\le G$ and $\chi\in \widehat{H}$.  Then $\vert \widehat{G}_{\chi} \vert = \lbrack G : H \rbrack$.
\end{cor}

\begin{proof}
Let $x_{1},x_{2},\ldots,x_{n}\in G$ such that $x_{1}\notin H$ and $x_{k}\notin H_{k-1}$, $1\le k \le n$, where $H_{0}=H$, $H_{n}=G$, and $H_{k}$ is the subgroup generated by $H_{k-1}$ and $x_{k}$.  Such elements exist because $G$ is a finite group.

The result will be demonstrated by induction.  The proof of Proposition \ref{Hext} shows for $\chi \in \widehat{H}$ that $\vert (\widehat{H_{1}})_{\chi} \vert = \lbrack H_{1}: H\rbrack$.  Assume that $\vert (\widehat{H_{k}})_{\chi}\vert = \lbrack H_{k} : H \rbrack$, with $1\le k\le n-1$.  Appealing again to the proof of Proposition \ref{Hext} it is clear that for any $\psi \in (\widehat{H_{k}})_{\chi}$,
\begin{equation*}
\vert (\widehat{H_{k+1}})_{\psi} \vert = \lbrack H_{k+1}: H_{k} \rbrack.
\end{equation*}

\noindent
But each element of $(\widehat{H_{k+1}})_{\psi}$, $\psi \in (\widehat{H_{k}})_{\chi}$, is an element of $(\widehat{H_{k+1}})_{\chi}$, so
\begin{equation*}
\vert (\widehat{H_{k+1}})_{\chi} \vert = \lbrack H_{k+1}: H_{k} \rbrack \, \lbrack H_{k}: H \rbrack = \lbrack H_{k+1} : H \rbrack.
\end{equation*}

\noindent
The $k=n$ instance of the induction statement is the conclusion of the proposition.
\end{proof}

\begin{cor} \label{Gchisum}
Let $H\le G$ and $\chi\in \widehat{H}$.  Then,
\begin{equation*}
\sum_{\psi \in \widehat{G}_{\chi}} \psi(g) = \begin{cases} \lbrack G:H \rbrack \chi(g), & g\in H, \\ 0, & \text{otherwise.} \end{cases}
\end{equation*}
\end{cor}

\begin{proof}
By the definition of $\widehat{G}_{\chi}$, $\psi(g)=\chi(g)$ whenever $g\in H$.  Hence, the claimed formula for $g\in H$ follows from Corollary \ref{Gchi}.  It is, therefore, sufficient to prove that if $g\notin H$, then $\sum_{\psi \in \widehat{G}_{\chi}} \psi(g) = 0$. Let $\lbrace x_{k} \rbrace_{k=1}^{n}$ and $H_{k}$, $0\le k\le n$, be as in the proof of Corollary \ref{Gchi}.  Recall from the proof of Proposition \ref{Hext} that each extension of $\varphi \in \widehat{H_{k}}$ to an element $\psi_{j}$ of $H_{k+1}$ is defined by
\begin{equation*}
\psi_{j}(x_{k+1}^{\ell}h) = \alpha_{j}^{\ell} \varphi(h),
\end{equation*}

\noindent
where $h\in H_{k}$ and, letting $m=\lbrack H_{k+1}:H_{k}\rbrack$, $\alpha_{j}$ is an $m$th root of $\varphi(x_{k+1}^{m})$. (Note that $x_{k+1}^{m}\in H_{k}$.)  It follows that
\begin{equation*}
\sum_{\psi \in (\widehat{H_{k+1}})_{\varphi}} \psi(g) = \varphi(h) \sum_{j=0}^{m-1} \alpha_{j}^{\ell},
\end{equation*}

\noindent
where $g = x^{\ell} h$, for some $h\in H_{k}$ and $0\le \ell < m$.  If $\ell\neq 0$, i.e., $g\notin H_{k}$, the sum is zero because $\lbrace \alpha_{j}\rbrace_{0}^{m-1}$ consists of a constant multiple of the $m$th roots of unity.  The result follows by an induction argument similar to that of Corollary \ref{Gchi}.
\end{proof}

The preceding corollaries set the stage for descriptions of sampling and upsampling over a subgroup $H\le G$ in terms of the DFT.

\begin{proposition} \label{DFTdownup}
Let $G$ be a finite abelian group with subgroup $H$.  Then
\renewcommand{\labelenumi}{(\roman{enumi})}
\begin{enumerate}
\item For $f\in \ell(H)$,
\begin{equation*}
\widehat{\mathcal{S}^{*}_{H} f}(\chi) = \hat{f}(\chi\vert_{H}), \quad \chi \in \widehat{G}.
\end{equation*}

\item For $f\in \ell(G)$,
\begin{equation*}
\widehat{\mathcal{S}_{H} f}(\chi) = \frac{1}{\lbrack G: H\rbrack} \sum_{\psi \in \widehat{G}_{\chi}} \hat{f}(\psi), \quad \chi\in \widehat{H}.
\end{equation*}
\end{enumerate}
\end{proposition}

\begin{proof}
Let $f\in \ell(H)$ and observe that
\begin{equation*}
\widehat{\mathcal{S}^{*}_{H} f}(\chi) = \sum_{g\in G} (\mathcal{S}^{*}_{H} f)(g) \overline{\chi(g)} = \sum_{h\in H} f(h) \overline{\chi(h)} = \hat{f}(\chi\vert_{H}).
\end{equation*}

\noindent
Given $f\in \ell(G)$, then
\begin{align*}
\widehat{\mathcal{S}_{H} f}(\chi) &= \sum_{h\in H} f(h) \overline{\chi(h)} \\
&= \sum_{g\in G} f(g) \frac{1}{\lbrack G:H\rbrack} \sum_{\psi \in \widehat{G}_{\chi}} \overline{\psi(g)} \\
&= \frac{1}{\lbrack G:H\rbrack} \sum_{\psi \in \widehat{G}_{\chi}} \hat{f}(\psi),
\end{align*}

\noindent
where the second equality follows from Corollary \ref{Gchisum}.
\end{proof}

\section{Convolutional systems for $\ell(G)$}

The results of the previous section make it possible to apply the approach of \cite{FJKO} for the study of convolutional systems in $\ell(\Bz/d\Bz)$ to convolutional systems in $\ell(G)$, where $G$ is an arbitrary finite abelian group.  The major tools required are the modulated filter representation (which rests on the sampling formulas of Proposition \ref{DFTdownup}) and Theorem \ref{FJKOtheorem12}.  The formal definition of a convolutional system for $\ell(G)$ follows.

\begin{definition}
Let $\lbrace f_{m} \rbrace_{m=0}^{n-1} \subset \ell(G)$, where $G$ is a finite abelian group.  Given a subgroup $H$ of $G$, the collection
\begin{equation*}
X_{H}\left ( \lbrace f_{m} \rbrace_{m=0}^{n-1} \right ) = \lbrace T_{h} f_{m} : h\in H, \; 0\le m \le n-1 \rbrace
\end{equation*}

\noindent
will be referred to as the \emph{convolutional system generated by $\lbrace f_{m}\rbrace_{m=0}^{n-1}$ with sampling over $H$}.
\end{definition}

\noindent
Note that the frame operator of such a system coincides with the filterbank frame operator described in Section \ref{section2} and depicted in Figure \ref{filterbank}.

\subsection{The modulated filter representation}

The purpose of the modulated filter representation is to generate a factorization of the filterbank synthesis operator, $L^{*}$, which will lead to a new system with equivalent frame properties whose synthesis operator possesses a block diagonal representation.

Towards this end, consider $y=y_{0} \oplus \cdots \oplus y_{n-1} \in \bigoplus_{m=0}^{n-1} \ell(H)$ under the action of $L^{*}$.  Applying $\mathcal{F}_{G}$ to $L^{*}y$ yields
\begin{equation*}
(\mathcal{F}_{G} L^{*} y ) (\psi) = \sum_{m=0}^{n-1} \hat{f}_{m} (\psi) \, \hat{y}_{m} (\psi \vert_{H}), \quad \psi \in \widehat{G},
\end{equation*}

\noindent
by Proposition \ref{DFTdownup}.  Alternatively, given $\chi\in \widehat{H}$,
\begin{equation} \label{lstareq}
(\mathcal{F}_{G} L^{*} y ) (\psi) = \sum_{m=0}^{n-1} \hat{f}_{m} (\psi) \, \hat{y}_{m} (\chi), \quad \psi \in \widehat{G}_{\chi}.
\end{equation}

\noindent
Notice that the value of $\mathcal{F}_{G} L^{*} y$ for each $\psi \in \widehat{G}_{\chi}$, $\chi \in \widehat{H}$, depends only on the Fourier transform of the components of $y$ at $\chi$.  There is no implicit ordering for the elements of $\widehat{G}$ and $\widehat{H}$, but, for the sake of a more explicit matrix represenation of the filterbank analysis operator, let the elements of $\widehat{H}$ be enumerated as $\widehat{H}=\lbrace \chi_{\ell} : 1\le \ell \le \vert H\vert \rbrace$ and those of $\widehat{G}$ in terms of $\widehat{G}_{\chi_{\ell}} = \lbrace \psi_{\ell,j} : 1\le j \le \lbrack G:H\rbrack \rbrace$.  (Notice that the latter enumeration is justified by Corollary \ref{Gchi}.)  Under these conventions, \eqref{lstareq} may be formulated as the matrix product
\begin{equation} \label{Hmodeq}
\frac{1}{\sqrt{\vert G\vert}} \begin{bmatrix} \mathcal{F}_{G} L^{*}y (\psi_{\ell,1}) \\ \vdots \\ \mathcal{F}_{G} L^{*}y (\psi_{\ell,\lbrack G:H\rbrack}) \end{bmatrix} = H^{*}_{\mathrm{mod}} (\ell) \, \frac{1}{\sqrt{\vert H\vert}} \begin{bmatrix} \hat{y}_{0}(\chi_{\ell}) \\ \vdots \\ \hat{y}_{n-1}(\chi_{\ell}) \end{bmatrix}, \quad 1\le \ell \le \vert H\vert,
\end{equation}

\noindent
where $H^{*}_{\mathrm{mod}}(\ell)$ is a $\lbrack G:H\rbrack \times n$ matrix whose $(j, m)$ entry is given by
\begin{equation*}
\lbrack H^{*}_{\mathrm{mod}} (\ell) \rbrack (j, m) = \frac{1}{\sqrt{\lbrack G:H\rbrack}} \hat{f}_{m}(\psi_{\ell,j}).
\end{equation*}

\noindent
The additional constants present in \eqref{Hmodeq} account for the fact that $\mathcal{F}_{G}$ and $\mathcal{F}_{H}$, as defined in Section \ref{section3}, are not unitary.  Collecting the $\vert H \vert$ equations represented by \eqref{Hmodeq} into a single matrix equation reveals the block diagonal nature of the modulated filter representation,
\begin{footnotesize}
\begin{equation*}
\frac{1}{\sqrt{\vert G\vert}}\!\!
\begin{bmatrix} \begin{array}{c}
\mathcal{F}_{G} L^{*}y (\varphi_{1,1}) \\ \vdots \\ \mathcal{F}_{G} L^{*}y (\varphi_{1,\lbrack G:H\rbrack \rbrace}) \\ \hline \mathcal{F}_{G} L^{*}y (\varphi_{2,1}) \\
\vdots\\ \mathcal{F}_{G} L^{*}y (\varphi_{2,\lbrack G:H\rbrack \rbrace}) \\ \hline
\\ \vdots \\ \\ \hline \mathcal{F}_{G} L^{*}y (\varphi_{\vert H\vert,1}) \\ \vdots \\ \mathcal{F}_{G} L^{*}y (\varphi_{\vert H\vert,\lbrack G:H\rbrack \rbrace}) \\
\end{array}\end{bmatrix}
\!\!=\!\! \frac{1}{\sqrt{\vert H\vert}}\!\!
\begin{bmatrix}\begin{array}{cccc}
H^{*}_{\mathrm{mod}}(\chi_{1})&\dots&0 \\ \vdots&\ddots&\vdots \\ 0&\dots& H^{*}_{\mathrm{mod}}(\chi_{\vert H\vert})
\end{array}\end{bmatrix}
\begin{bmatrix} \begin{array}{c}
\mathcal{F}_{H} y_0(\chi_{1}) \\ \vdots \\ \mathcal{F}_{H} y_{n\!-\!1}(\chi_{1}) \\ \hline \mathcal{F}_{H} y_0(\chi_{2}) \\ \vdots \\ \mathcal{F}_{H} y_{n\!-\!1}(\chi_{2}) \\ \hline \\ \vdots \\ \\ \hline \mathcal{F}_{H}y_0(\chi_{\vert H\vert}) \\ \vdots \\ \mathcal{F}_{H} y_{n\!-\!1}(\chi_{\vert H\vert})
\end{array}\end{bmatrix}.
\end{equation*}
\end{footnotesize}

\noindent
The block matrix on the right-hand side of the last equality is called the \emph{block adjoint modulated filter matrix} and will be denoted by $H_{\mathrm{mod}}^{*}$.  Notice that $U_{1}=\vert G\vert^{-1/2} \mathcal{F}_{G} \in \mathcal{U}(\ell(G))$ and $U_{2}=\vert H\vert^{-1/2} \mathcal{F}_{H} \in \mathcal{U}(\ell(H))$, which implies
\begin{equation} \label{unitaryeq}
L^{*} = U_{1}^{-1} H_{\mathrm{mod}}^{*} (\oplus_{\ell} U_{2}).
\end{equation}

\noindent
It should be noted that given another ordering of the elements comprising $\widehat{G}$ or  $\widehat{H}$ one could introduce additional permutation matrices to achieve the above block diagonalization.  Equation \eqref{unitaryeq} motivates the definition of the $\lbrack G:H\rbrack$ collections,
\begin{equation*}
Y_{\ell} = \lbrace y_{m,\ell} \rbrace_{m=0}^{n-1} \subset \ell(\Bz/\lbrack G:H\rbrack \Bz), \quad 1\le \ell \le \vert H\vert,
\end{equation*}

\noindent
where $y_{m,\ell}(j) = \lbrack G:H\rbrack^{-1/2} \hat{f}_{m}(\psi_{\ell,j})$.  As a result of this definition, $H^{*}_{\mathrm{mod}}(\ell)$ is the synthesis operator corresponding to $Y_{\ell}$.  Moreover, since $L^{*}L$ and $H_{\mathrm{mod}}^{*} H_{\mathrm{mod}}$ are unitarily equivalent, it follows that:
\begin{enumerate}
\item the frame bounds for $X_{H}(\lbrace f_{m}\rbrace_{m=0}^{n-1})$ are the minimum of the lower frame bounds and the maximum of the frame bounds for the collections $Y_{\ell}$;

\item the frame potential of $X_{H}(\lbrace f_{m}\rbrace_{m=0}^{n-1})$ is equal to the sum of the frame potential of the collections $Y_{\ell}$, $1\le \ell \le \vert H\vert$.
\end{enumerate}

\noindent
The arguments used to justify these facts can be found in \cite{FJKO} and, therefore, will not be repeated here.

\subsection{The main result}

\begin{theorem} \label{mainthm}
Let $\lbrace a_{m} \rbrace_{m=0}^{n-1} \subset \Br$ be such that $a_{0} \ge a_{1} \ge \cdots \ge a_{n-1}>0$.  Let $G$ be a finite abelian group and $H$ a subgroup of $G$ with $n\ge \lbrack G:H\rbrack$.  Denote by $m_{0}$ the smallest index $0\le m \le \lbrack G:H\rbrack -1$ such that
\begin{equation}  \label{FFineqH}
(\lbrack G:H\rbrack -m)a_{m}^{2} \le \sum_{j=m}^{n-1} a_{j}^{2}.
\end{equation}

\noindent
If $X_{H}( \lbrace f_{m} \rbrace_{m=0}^{n-1}) \subset \ell(G)$ is a local minimizer of the frame potential under the constraint that $\Vert f_{m}\Vert = a_{m}$, $0\le m \le n-1$, then $X_{H}(\lbrace f_{m}\rbrace_{m=0}^{n-1})$ may be divided into two mutually orthogonal subcollections of $\ell(G)$: $X_{H}(\lbrace f_{m}\rbrace_{m=0}^{m_{0}-1})$, which consists of mutually orthogonal, nonzero vectors, and $X_{H}(\lbrace x_{m}\rbrace_{m=m_{0}}^{n-1})$, which is a tight frame for its $\vert H\vert (\lbrack G:H\rbrack -m_{0})$-dimensional span.  In particular, if $m_{0}=0$ then $X_{H}(\lbrace f_{m}\rbrace_{m=0}^{n-1})$ is a tight frame for $\ell(G)$.
\end{theorem}

\begin{proof}
Consider the collections $Y_{\ell}$, $1\le \ell \le \vert H\vert$, defined in the previous section and observe that, up to a constant multiple, the Fourier coefficients of a given filter $f_{m}$ are distributed as the values of the elements $y_{m,\ell}(j)$, $1\le \ell \le \vert H\vert$, $1\le j\le \lbrack G:H\rbrack$.  Hence, the constraint $\Vert f_{m}\Vert = a_{m}$ can be restated as
\begin{equation*}
\sum_{\ell=1}^{\vert H\vert} \Vert y_{m,\ell} \Vert^{2} = \vert H\vert a_{m}^{2}, \quad 0\le m \le n-1.
\end{equation*}

\noindent
Referring to Theorem \ref{FJKOtheorem12} with $N=\lbrack G:H\rbrack$ it is not difficult to see that if the hypotheses of Theorem \ref{mainthm} are satisfied, then those of Theorem \ref{FJKOtheorem12} will be satisfied with the same $m_{0}$.  Theorem \ref{FJKOtheorem12} thus provides a decomposition of each $Y_{\ell}$, $1\le \ell \le \vert H\vert$, into mutually orthogonal subcollections $\lbrace y_{m,\ell} \rbrace_{m=0}^{m_{0}-1}$ and $\lbrace y_{m,\ell} \rbrace_{m=m_{0}}^{n-1}$, the former collection consisting of mutually orthogonal vectors and the latter comprising a tight frame for its span.  In particular, for each $1\le \ell \le \vert H\vert$, $\Vert y_{m,\ell}\Vert^{2}=a_{m}^{2}$ for $m<m_{0}$ and
\begin{equation} \label{TFbound}
\sum_{m=m_{0}}^{n-1} \Vert y_{m,\ell}\Vert^{2} = \sum_{m=m_{0}}^{n-1} a_{m}^{2}, \quad m_{0} \le m \le n-1.
\end{equation}

\noindent
Let $F_{Y_{\ell}}$ denote the frame operator of $Y_{\ell}$ on $\ell(\Bz/\lbrack G:H\rbrack \Bz)$ and $F_{Y} = H_{\mathrm{mod}}^{*} H_{\mathrm{mod}}$ the frame operator induced by $H_{\mathrm{mod}}^{*}$ on $\ell(G)$.  It follows from the above decomposition that each $y_{m,\ell}$ is an eigenvector of $F_{Y_{\ell}}$ and that the eigenvalue is independent of $\ell$.  If $m<m_{0}$, then the eigenvalue is $a_{m}^{2}$, while if $m\ge m_{0}$, the eigenvalue is given by \eqref{TFbound}.  As a result, $\oplus_{\ell} y_{m,\ell}$ is an eigenvector of $F_{Y}$ with the corresponding eigenvalue.  By construction, $\oplus_{\ell} y_{m,\ell}$ is equal to $ \sqrt{\vert H\vert} U_{1} f_{m}$ and
\begin{equation*}
F_{Y} = H^{*}_{\mathrm{mod}} H_{\mathrm{mod}} = U_{1} L^{*} L U_{1}^{-1} = U_{1} F U_{1}^{-1},
\end{equation*}

\noindent
where $F$ is the frame operator of the collection $X_{H}$.  Since $U_{1} f_{m}$ is an eigenvector of $F_{Y}$, it follows that $f_{m}$ is an eigenvector of $F$.  Moreover, the eigenvalue of $f_{m}$ for $F$ must equal that of $\oplus_{\ell} y_{m,\ell}$ for $F_{Y}$, while the symmetry of the collection $X_{H}$ forces each translate $T_{h}f_{m}$, $h\in H$, to be an eigenvector with the same eigenvalue.  Finally, note that the claimed orthogonality of the subcollections $X_{H}(\lbrace f_{m}\rbrace_{m=0}^{m_{0}-1})$ and $X_{H}(\lbrace x_{m}\rbrace_{m=m_{0}}^{n-1})$ (along with the mutual orthogonality of the vectors in the former subcollection) is inherited from the above decomposition of $Y_{\ell}$, since the frame operators $F_{Y}$ and $F$ are unitarily equivalent. This completes the proof.
\end{proof}

The final result of this section examines underdetermined systems and follows from an argument analogous to that used in \cite{FJKO}.

\begin{cor} \label{underdetermined}
Let $\lbrace a_{m} \rbrace_{m=0}^{n-1} \subset \Br$ be such that $a_{0} \ge a_{1} \ge \cdots \ge a_{n-1}>0$.  Let $G$ be a finite abelian group and $H$ a subgroup of $G$ with $n\le \lbrack G:H\rbrack$.  If $X_{H}(\lbrace f_{m} \rbrace_{m=0}^{n-1}) \subset \ell(G)$ is a local minimizer of the frame potential under the constraint that $\Vert f_{m}\Vert = a_{m}$, $0\le m \le n-1$, then $X_{H}(\lbrace f_{m}\rbrace_{m=0}^{n-1})$ is an orthogonal sequence in $\ell(G)$.
\end{cor}


\section{Directions for further study}

Recall that one consequence of the characterization results of Theorems \ref{FFthm} and \ref{mainthm} is the existence of tight frames with a certain structure.  In the case of Theorem \ref{FFthm} the structure imposed is limited to the norms of the frame elements, while in Theorem \ref{mainthm} the frames studied are required to consist of the translates of elements of prescribed norms over a subgroup.  In each case, the fundamental frame inequality describes when the imposed structure permits tight frames.  It is natural to wonder what other structures lead to similar descriptions of tight frames and whether some version of the fundamental frame inequality will always appear in the characterization.  The following specific questions are motivated by this idea.
\begin{enumerate}
\item \emph{Is it possible to extend the characterization of tight frames in terms of the frame potential to convolutional systems for $\ell(G)$, where $G$ is an arbitrary finite group?}

\item \emph{What other symmetries or structures of systems in a finite dimensional Hilbert space lead to similar characterizations of tight frames in terms of the frame potential?}
\end{enumerate}

In a certain sense, the frame potential measures the orthogonality of a given collection of vectors and, as such, minimizers of the frame potential can be thought to represent maximally orthogonal collections.  Recent interest in equiangular frames seems to take a step away from this concept; however, it is reasonable to ask whether potential methods can be adapted sufficiently to produce existence results for equiangular frames.

\section*{Acknowledgments}

The first author would like to thank Paul Koester for insightful exchanges on the Fourier analysis of finite groups.

\bibliographystyle{amsalpha}

\end{document}